\documentclass[twoside,12pt]{article}
\usepackage{amsmath,amssymb,epsfig, amsthm}
\date{}
\newtheorem{Lemma}{LEMMA}[section]
\newtheorem{Corollary}[Lemma]{COROLLARY}
\newtheorem{Theorem}[Lemma]{THEOREM}
\newtheorem{Proposition}[Lemma]{PROPOSITION}

\newcommand{\bnum}{\begin{enumerate}}
\newcommand{\enum}{\end{enumerate}}
\newcommand{\bi}{\begin{itemize}}
\newcommand{\ei}{\end{itemize}}
\newcommand{\btab}{\begin{tabular}}
\newcommand{\etab}{\end{tabular}}
\newcommand{\beq}{\begin{eqnarray*}}
\newcommand{\eeq}{\end{eqnarray*}}
\newcommand{\beqn}{\begin{eqnarray}}
\newcommand{\eeqn}{\end{eqnarray}}

\newcommand{\bq}{\begin{equation}}
\newcommand{\eq}{\end{equation}}
\newcommand{\CA}{{\mathcal A}}

\newcommand{\CK}{{\mathcal K}}
\newcommand{\CL}{{\mathcal L}}

\def\phi{\varphi}
\def\epsilon{\varepsilon}

\newcommand{\BR}{\mathbb R}

\newcommand{\kasten}{\vbox{\hrule height 8pt width 8.6pt depth -7.4pt
    \hbox{\vrule width 0.6pt height 7.4pt
    \kern 7.4pt \vrule width 0.6pt height 7.4pt}
    \hrule height 0.6pt width 8.6pt}}
\newcommand{\ok}{\hfill\kasten}
\newcommand{\bpf}{\begin{Proof}}
\newcommand{\epf}{\ok\end{Proof}\bigskip\noindent}
\newcommand{\bthm}{\begin{Theorem}}
\newcommand{\ethm}{\end{Theorem}}
\newcommand{\ble}{\begin{Lemma}}
\newcommand{\ele}{\end{Lemma}}
\newcommand{\bprop}{\begin{Proposition}}
\newcommand{\eprop}{\end{Proposition}}
\newcommand{\bcor}{\begin{Corollary}}
\newcommand{\ecor}{\end{Corollary}}

\begin{document}
\title{Semiaffine stable planes}
\author{Rainer L\"owen and Markus J. Stroppel}
\maketitle
\thispagestyle{empty}

\begin{abstract}

A locally compact stable plane of positive topological dimension 
will be called semiaffine if for every line $L$ 
and every point  $p$ not in $L$ there is at most one line passing through $p$ and disjoint 
from $L$. We show that then the plane is either an affine or projective plane or a 
punctured projective plane (i.e., a projective plane with one point deleted). 
We also compare this with the situation in general linear spaces (without topology), 
where P. Dembowski showed that the analogue of our main result is true for finite 
spaces but fails in general.
\\


\end{abstract}

\section{Introduction}\label{intro}

We recall the \it basic notions, \rm compare \cite{diss}, 
\cite{central}, \cite{handb},  \cite{adv}. 
A \it stable plane \rm is a pair $(M,\CL)$ consisting of a set $M$ of \it points \rm 
and a set $\CL$ of subsets $L\subseteq M$, called \it lines\rm, such that any two points 
$p,q \in M$ are joined by a unique line $L= p \vee q \in \CL$. Moreover, it is required 
that both $M$ and $\CL$ carry locally compact topologies of positive topological dimension 
such that the operation $\vee \colon M \times M \to \CL$ of \it joining \rm is continuous 
and the opposite operation $\wedge$ of \it intersection \rm sending a pair of intersecting 
lines to their point of intersection is continuous and its domain of definition is open.
The last mentioned property is called \it stability of intersection\rm . It distinguishes
stable planes from spatial geometries, where pairs of intersecting lines may be 
approximated by pairs of skew lines. To avoid trivialities,
it is also assumed that there is a \it quadrangle\rm, i.e., four points no three of 
which are on one line.

If two lines always intersect, then we have a \it topological projective plane\rm .  
This happens
if and only if $M$ is compact. See \cite{CPP} for many examples. An abundance 
of non-projective examples is obtained by taking the trace geometry induced on any proper open
subset of the point set in a topological projective plane. Apart from these 
\it (projectively)  embeddable \rm 
examples, there are non-embeddable ones. They are generally harder to construct, and 
this explains why only two-dimensional examples are known, compare \ref{example} below. 

The \it pencil \rm of all lines passing through a point $p$ is denoted $\CL_p$; it 
is always compact and connected, see \cite{diss}, 1.17, 1.14. 
As a consequence, a non-compact line $L$ is disjoint 
from at least one line through any point $p \notin L$. This is seen by looking at the map 
$x \to x \vee p$ from $L$ to $\CL_p$. 

In \cite{central}, we studied special cases of Euclid's \it parallel axiom\rm, 
which stipulates that given a line $L$ and a point $p \notin L$, there exists exactly one line 
$K \in \CL_p$ that is disjoint from $L$. This line is then said to be \it parallel \rm to $L$.
Also every line is considered as being parallel to itself.
If the parallel axiom holds universally, then the plane is called an \it affine plane. \rm
If the parallel axiom is satisfied for a fixed line $A$ and for each point not on $A$, 
then we say that the line $A$ is \it affine. \rm 

An affine line all of whose parallels 
are also affine will be called a \it biaffine \rm line. Not all affine lines are 
biaffine. Indeed, if we take a projective plane $(P,\CL)$ and delete from $P$ a closed 
subset of some line $W$, then what we obtain is called an 
\it almost projective plane\rm. If we have 
deleted a proper subset $X \subseteq W$ containing more than one point, 
then any line $A\ne W$ originally containing a point of $X$ will be affine, 
and the remainder $W\setminus X$ is parallel to $A$, but is not affine.
So $A$ is affine, but not biaffine. Two intersecting affine lines may 
have a common parallel, as shown by lines $A_1$ and $A_2$ of the kind just considered, 
the common parallel being $W\setminus X$. However, we observe:

\ble \label{equiv}
On the set $\CA$ of all affine lines of a stable plane, the parallelity relation 
is an equivalence.
\ele

\bpf 
We have to show that two parallels of an affine line  do 
not intersect. Indeed, this is part of  the definition of an affine line.
\epf

\ble\label{cts}
There is a continuous map $P \colon \CA \times M \to \CL$ that sends a pair $(A,p)$ to the 
unique parallel of $A$ passing through $p$. 
\ele

\bpf
The spaces $M$ and $\CL$ are second countable by \cite{diss}, 1.9, hence it suffices to 
show sequential continuity. So let $(A_n,p_n) \to (A,p)$ in $\CA \times M$. 
If $P(A_n,p_n)$ converges to a line $L$, then stability of intersection implies that 
either $p\in A$ and $L = A$ or $p\notin A$ and 
$L$ is a line passing through $p$ which is disjoint from $A$. 
So in either case, $L = P(A,p)$. Now by \cite{diss}, 1.17,
the set of all lines containing one of the points $p_n$ or $p$ is compact, hence
every subsequence of $P(A_n,p_n)$ has a subsequence converging to $P(A,p)$. 
This implies convergence of the entire sequence.
\epf

A line is $L$ called \it projective \rm if it meets every other line. 
A line is projective if and only if it is compact, because we have a continuous and 
open embedding $p \mapsto p \vee q$ of $L$ into the compact and connected pencil $\CL_q$ of any 
point $q \notin L$.

Finally, we need the notion of a \it pointwise coaffine line\rm, 
introduced in \cite{central}.
These can be characterized as non-projective lines intersecting only all projective lines. 
We refer to \cite{central} for an explanation of the term `pointwise coaffine', which 
describes how these lines appear in the so-called opposite plane, a kind of weak dual. 
In a \it punctured projective plane\rm, i.e., a projective plane with just one 
point deleted, all lines are either projective or pointwise coaffine, and both types occur. 
In fact, this property characterizes punctured projective planes as is easily seen.
We give a proof, because this fact does not seem to be recorded in the literature:
 
\bprop\label{coaffine}
If every line of a stable plane is either projective or pointwise coaffine then the 
plane is either a projective or a punctured projective plane. The converse is also true.
\eprop

\bpf We may assume that non-projective (hence non-compact) lines exist. 
Every non-compact line has at least one parallel through any given point outside.
That line is not projective and must be pointwise coaffine, hence it is 
the unique non-compact line passing through that point. 
Thus we have a set of pairwise disjoint biaffine lines covering the point set. 
Adding a point at infinity, we obtain a topological projective plane, compare
\cite {central}, Theorem 2.2. 
\epf

In the terminology of \cite{central},  
punctured projective planes are also called \it coaffine planes\rm.
Here is another related fact, taken from \cite{central}, 1.3.

\bprop\label{pencil}
If some pencil $\CL_p$ of a stable plane consists entirely of projective lines, 
then the plane is projective.
\eprop

\bpf
Every line $L$ not containing $p$ is homeomorphic to the compact pencil $\CL_p$ via join 
and intersection, hence $L$ is projective.
\epf

Our focus here will be on \it semiaffine \rm stable planes, that is, planes where for 
every non-incident point-line pair $(p,L)$ there is at most one line containing 
$p$ and disjoint from $L$. We have the following

\bprop\label{weakly aff}
In a semiaffine stable plane, every line is either affine or projective.
\eprop

\bpf
If all lines passing through $p \notin L$ intersect $L$, then $L$ is homeomorphic to 
the compact pencil $\CL_p$, hence $L$ is projective. If this never happens for a given line $L$,
then $L$ is affine.
\epf

\section{Semiaffine planes}\label{saff}

\ble\label{affine} 
If a projectively embeddable stable plane $(M,\CL)$ contains two intersecting 
affine lines, one of which is biaffine, then $(M,\CL)$ is an affine plane. \ele

\bpf
Suppose that $(M,\CL)$ is embedded in the topological projective plane $(P,\CK)$. 
Recall from the introduction that this means that $M$ is an open subset of $P$, 
and that the lines $L \in \CL$
are the nonempty intersections $L = K \cap M$, $K \in \CK$. Note that then $L$ 
contains more than one point.

For $L \in \CL$, let 
$\bar L \in \CK$ be the unique line containing $L$. If $\bar L \setminus L$ contains 
distinct points $a,b$, then any point $p\in M \setminus L$ lies on two distinct 
lines $p \vee a$ and $p \vee b$ disjoint from $L$. Hence $L$ is affine if and only if 
$\bar L \setminus L$ consists of a single point $\infty_L$. The parallels of $L$ 
in the stable plane $(M,\CL)$ are then 
precisely the lines $K \in \CL$ such that $\bar K$ contains $\infty_L$. If $L$ 
is biaffine, it follows that $\bar K \setminus K = \{\infty_L\}$ holds for these lines. 
Let $\CA_L$ be the set of all parallels $K$ of $L$. Passing to the corresponding 
lines $\bar K$, we obtain a subset $\overline {\CA_L} \subseteq \CK_{\infty_L}$, and 
then we see that $M$ is the union of $\overline {\CA_L}$ minus the point $\infty_L$. 

Now let $B\notin \CA_L$ be another affine line. Then $\infty_B \ne \infty_L$, and the 
line $W\in \CK$ joining these two points does not belong to $\overline {\CA_L}$, 
because it contains two points outside $M$. On the other hand, all other lines $C$ in the pencil
of $\infty_L$ must belong to $\overline {\CA_L}$ because their intersection $q$ with 
$\bar B$ belongs to $B \subseteq M$, whence $q\vee \infty_L$ induces a parallel of $L$. 
It follows that $\overline {\CA_L} = \CK_{\infty_L} \setminus \{W\}$,
and hence that $M = P \setminus W$ is an affine plane. In fact, $(P,\CK)$ is the 
projective completion of this affine plane.
\epf

In \cite{central}, 5.4, a very special example of a stable plane is constructed. 
It has point set $M = \BR^2$ and three types of lines, all homeomorphic to $\BR$: 
one special line $C$ (the $y$-axis), the lines not meeting $C$, and the lines meeting $C$.
Precisely the lines of the third type are affine. They are in fact biaffine, because 
their parallels are obtained by translation in $y$-direction. The plane is 
not affine, because the lines of the second type are not affine. Using Lemma \ref{affine},
we infer the following

\bcor\label{example} The two-dimensional stable plane referred to above is an example of a 
non-embeddable plane.
\ok
\ecor

The plane mentioned can be extended by adding a point at infinity to every line of the third 
type, see \cite{central}, 5.4. The added points together form an additional line $D$.
One obtains an example of a stable plane containing two pointwise 
coaffine lines $C$ and $D$ without being a coaffine (i.e., punctured projective) plane. This 
property shows directly that the extended plane is non-embeddable. We could have 
told this before, since an open embedding of the extended plane would induce an open 
embedding of the original one.

\ble\label{intersect}
Suppose that $A$ and $B$ are two intersecting biaffine lines in a stable plane. Then 
every parallel of $A$ intersects every parallel of $B$.
\ele

\bpf
Let $A'$ and $B'$ be parallels of $A$ and $B$, respectively. Then
$A'$ and $B'$ are affine by assumption. If they are disjoint, i.e., parallel,
then this violates Lemma \ref{equiv}.
\epf

This is false if $A$ and $B$ are merely affine. Indeed, if then $A'$ 
intersects $B'$, deleting the intersection point from the given plane yields a counterexample.

\bprop
If $A$ and $B$ are two intersecting biaffine lines in a stable plane $(M,\CL)$, then the 
map $A\times B \to M$ that sends $(a,b)$ to the point $P(a,B) \wedge P(b,A)$ is a homeomorphism.
\eprop

\bpf
By Lemma \ref {intersect} we have an  inverse map sending 
$p \in M$ to $(P(B,p)\wedge A, P(A,p) \wedge B)$; it is continuous 
by Lemma \ref{cts}.
\epf

\bcor \label{noproj}
If a stable plane contains two intersecting biaffine lines $A$, $B$, then it does not 
contain a projective line.
\ecor

\bpf
If $L$ is a projective line, then by Lemma \ref{intersect}
there is a continuous map $L \to B$ sending $p \in L$ to
$P(A,p) \wedge B$. This map is surjective because $L$ meets every line. 
However, $L$ is compact and $B$ is not. 
\epf

In none of the two preceding assertions, the assumption can be weakened 
from `biaffine' to `affine';
every almost projective plane that is neither projective nor punctured projective 
nor affine yields counterexamples. The following lemma is obvious from the definitions.

\ble \label{bi}
If all lines of a stable plane are either affine or projective, 
then every affine line is biaffine.
\ok
\ele

\bthm\label{main}
The semiaffine locally compact positive-dimensional stable planes are 
precisely the following:
\begin{enumerate}
\item affine planes
\item projective planes
\item punctured projective planes.
\end{enumerate}
\ethm

\bpf
Clearly, the planes listed are all semiaffine. For the converse assertion,
recall that all lines are either affine or projective by Proposition \ref{weakly aff}.
We only need to consider the mixed case, where both affine and projective lines exist. 
We claim that the plane is punctured projective. By Lemma \ref{bi}, all affine 
lines are biaffine, and by Lemma \ref{noproj}, no two of them intersect. We may then embed 
our plane in a compact projective plane by adding a point $\infty$ to the point set 
and replacing every affine line $A$ by $A\cup\{\infty\}$. See \cite{central}, Theorem 2.2 for 
details on the continuity properties of the extended plane. 

Instead of this completion argument,
we may use the notion of a coaffine point introduced in \cite{central}, as follows: 
By the previous remarks, every point is incident with precisely one affine line, and all 
other lines passing through that point are projective. This means that every point is 
coaffine. Then the plane is pointwise coaffine and hence coaffine, i.e., punctured projective. 
\epf

We take this opportunity to give a complete statement of Proposition 1.8 of \cite{central}, 
which was mutilated by the publisher after proofreading. 

\bprop
{\rm \cite{central}} Let $(M,\CL)$ be a stable plane whose lines are manifolds.

\item{(a)} If $C\in \CL$ is pointwise coaffine, then the complement $(M \setminus C, 
\CL\setminus \{C\})$ has a point set homeomorphic to Euclidean space $\BR^{2^{n+1}}$, 
where $0 \le n \le 3$, and all lines are homeomorphic to $\BR^{2^n}$.

\item {(b)}  If there is a second pointwise coaffine line $D$, then the opposite plane
$(M,\CL)^*$ has the same topological properties.
\eprop

The proof given in \cite{central} is correct. The opposite plane has point set $\CK$, 
the set of all compact lines of $(M,\CL)$, and its lines are the partial pencils $\CK_p$.

\section{Appendix: Linear spaces} 

All our assertions make sense in the more general situation of \it linear spaces\rm,
where there is no topology, and we shall examine which of them remain true. Some strong tools,
in particular those related to compactness, are no longer available, and examples are hard 
to construct. This is why at least one question remains open. Yet also some strong 
results are known.

We define a linear space to be a pair $(M,\CL)$ consisting of a set $M$ and a collection 
$\CL$ of subsets of cardinality at least 3, called lines, such that any two points are 
joined by a unique line, and such that every point is on at least 3 lines. 
Linear spaces include all projective or affine spaces, but if we assume the existence 
of affine or projective lines (defined in the same manner as before), then this 
indicates planar behaviour. Kreuzer \cite{kreuz} introduces a notion of semiaffine linear 
spaces that makes sense in higher definitions, but we shall not adopt his definition. 

A few results from Section \ref{saff} remain true without the topological assumptions. 
Notably, this holds for Lemma \ref{equiv} 
(parallelity is an equivalence among affine lines) and its consequence Lemma \ref{intersect}
(parallels of intersecting biaffine lines always intersect). The proofs are valid 
without change. All remaining results that we obtained for stable planes rely on 
topological arguments for their proofs, and we have stressed the places where 
this occurs in the previous sections. \\

One assertion that definitely fails is the analogue of Proposition \ref{weakly aff}, 
which asserts that a 
semiaffine linear space contains only affine and projective lines. An easy 
counterexample is the projectively embedded plane
obtained from any affine plane by adding a single point $x$ at infinity in 
the projective closure. 
Lines passing through $x$ are projective, but those not containing $x$ are not affine, 
because they do not possess parallels containing $x$. As a consequence, the analogue of 
Theorem \ref{main} fails. In addition to this standard counterexample, 
P. Dembowski \cite{dembo} constructed a plethora of other
counterexamples by a processs of free extension. 
Also Proposition \ref{pencil} fails without topology: in the above example, 
the pencil of $x$ consists 
of projective lines. We do not know the answer to the following\\

\bf PROBLEM: \rm Does Theorem \ref{main} hold for linear spaces with the 
stronger assumption that all 
lines are either affine or projective? \\

Things are somewhat different in finite linear spaces. There, 
existence of a projective or affine line 
implies that all pencils of points outside this line have equal cardinalities, and this 
may sometimes replace our compactness arguments. In fact, we have the following.

\bthm\label{finite}
Let $(M,\CL)$ be a finite linear space as defined above. If all lines are either 
affine or projective, then $(M,\CL)$ is an affine plane or a projective plane or a 
punctured projective plane.
\ethm

\bpf
Only the mixed case (with lines of both kinds) needs to be discussed. 
No line is affine and projective at the same time, and every affine line is biaffine. 
Let $q$ be the number of points on an affine line $A$. For every point $x$ outside $A$, 
we have the unique parallel to $A$ containing $x$. The other lines in the pencil $\CL_x$ 
meet $A$, hence the pencil has $q+1$ elements. Repeating this argument with another 
line parallel to $A$ one sees that every pencil has cardinality $q+1$. It follows 
that a line is affine or projective according as its cardinality is $q$ or $q+1$, respectively.

Now we can use the arguments from Lemmas \ref{intersect} and \ref{noproj} 
(with homeomorphisms replaced by bijections) to conclude that no projective lines 
exist if there are two intersecting affine lines. 
The only remaining possibility is that every point is on a unique affine line, and then by
adding a common point at infinity to these lines we obtain a projective plane. 
\epf

Theorem \ref{finite} is a special case of a result by N. Kuiper and P. Dembowski \cite{dembo}, 
which is much harder to prove. It asserts that the only finite semiaffine linear spaces are the
finite affine, projective or punctured projective planes and the finite affine planes 
extended by one point at infinity. In other words, the counterexample to Theorem 
\ref{main} exhibited above 
is the only one in the finite case.
A generalization allowing lines with only two points 
is given by J. Totten and P. de Witte \cite{totten}.

\bibliographystyle{plain}

\begin{thebibliography}{9}

\bibitem{dembo}
P. Dembowski,
Semiaffine Ebenen,
{\em Arch. Math.} 13, 120 -- 131, 1962.

\bibitem{handb}
T. Grundh\"ofer and R. L\"owen,
Linear topological geometries,
in: F. Buekenhout (ed.), 
{\em Handbook of Incidence Geometry}, Chapter 23, pp. 1255 -- 1324,
North Holland, Amsterdam, 1995.

\bibitem{kreuz}
A. Kreuzer,
Semiaffine spaces,
{\em J. Comb. Theory, Ser. A} 64, 63 -- 78, 1993.


\bibitem{diss}
R. L\"owen,
Vierdimensionale stabile Ebenen,
{\em Geometriae Dedicata} 5, 239 -- 294, 1976.


\bibitem{central}
R. L\"owen,
Central collineations and the parallel axiom in stable planes,
{\em Geometriae Dedicata} 10, 283 -- 315, 1981.

\bibitem{adv}
H. Salzmann,
Topological planes,
{\em Advances in Math.} 2, 1 -- 60, 1967.

\bibitem{CPP}
H. Salzmann, D. Betten, T. Grundh\"ofer, H. H\"ahl, R. L\"owen, and M. Stroppel,
{\em Compact projective planes}, 
Walter de Gruyter, Berlin, New York 1995.



\bibitem{totten}
J. Totten and P. de Witte,
On a Paschian condition for linear spaces,
{\em Math. Z.} 137, 173 -- 183, 1974.

\end{thebibliography}

\bigskip
\bigskip
\noindent{Rainer L\"owen\\ 
Institut f\"ur Analysis und Algebra\\
Technische Universit\"at Braunschweig\\
Universit\"atsplatz 2\\
38106 Braunschweig\\
Germany\\
r.loewen@tu-braunschweig.de\\ \\
Markus J. Stroppel\\
Fakult\"at 8\\
LExMath\\
Universit\"at Stuttgart\\
70550 Stuttgart \\
Germany\\
stroppel@mathematik.uni-stuttgart.de}

\end{document}